\documentstyle{article}
\begin{document}
\title{su$_q$(2)-Invariant  Harmonic Oscillator \footnote{presented at the $8th$ Colloquium "Quantum groups and Integrable
systems", Prague 17-19 June 1999.}}
\author{M. IRAC-ASTAUD \\
Laboratoire de Physique Th\'eorique de la Mati\`ere
Condens\'ee,\\
Universit\'e
Paris VII, 2, Place Jussieu,\\ F-75251 Paris Cedex 05, France\\
e-mail : mici@ccr.jussieu.fr
\and
C.\ QUESNE \thanks{Directeur de recherches FNRS}\\
Physique Nucl\'eaire Th\'eorique et Physique
Math\'ematique,\\
Universit\'e Libre
de Bruxelles, Campus de la Plaine CP229,\\ Boulevard du Triomphe, B-1050 Brussels,
Belgium \\
e-mail : cquesne@ulb.ac.be}

\date{}
\maketitle


\begin{abstract}
We propose a $q$-deformation of the su(2)-invariant Schr\"odinger equation of a
spinless particle in a central potential, which allows us not only to
determine a
deformed spectrum and the corresponding eigenstates, as in other approaches, but
also to calculate the expectation values of some physically-relevant
operators. Here
we consider the case of the isotropic harmonic oscillator and of the quadrupole
operator governing its interaction with an external field. We obtain the
spectrum
and wave functions both for $q \in $ R$^+$ and generic $q \in S^1$, and study the
effects of the $q$-value range and of the arbitrariness in the su$_q$(2) Casimir
operator choice. We then show that the quadrupole operator in $l=0$ states
provides
a good measure of the deformation influence on the wave functions and on the
Hilbert space spanned by them.
\end{abstract}

%
\section{Introduction}

Since the advent of quantum groups and quantum algebras, there has
been a lot of interest in deformations of the harmonic oscillator,
since the latter plays a central role in the investigation of many
physical systems. Various $q$-deformed versions of standard quantum
mechanics in the Schr\"odinger representation were proposed for the
oscillator by using either the ordinary differentiation operator
(see e.g.~\cite{macfarlane}), or a $q$-differentiation one (see
e.g.~\cite{ minahan}). Some works also involve non-commutative
objects (see e.g.~\cite{carow94}).
%

%
In previous papers~\cite{irac96,nous} and in this talk, we set up
an
  su$_q$(2)-invariant Schr\"odinger equation within the framework of
   quantum mechanics, using a representation of the
su$_q$(2) quantum algebra on the two-dimensional
sphere~\cite{rideau,irac98},
 which allows us not only to determine a deformed spectrum, as in other
 works,
 but also
  to calculate the expectation values of some physically-relevant operators.
%
%
In our approach, only the angular sector is deformed. This gives
rise to an appropriate change in the angular part of the scalar
product~\cite{irac98}, and to the substitution of the su$_q$(2)
Casimir operator for the su(2) one in the kinetic part of the
Hamiltonian while the potential part remains unchanged. The latter
step may be performed in various ways since there is no unique rule
for constructing the su$_q$(2) Casimir operator. Similarly, the
deforming parameter $q$ may be assumed either real and positive, or
on the unit circle in the complex plane (but different from a root
of unity), provided different scalar products are
used~\cite{irac98}. We will study the effects of these two choices
on the spectrum and on the value of the quadrupole moment.\par
%
%
In Section~\ref{sec:equation}, the Schr\"odinger equation of the
su$_q$(2)-invariant harmonic oscillator is introduced and solved.
In Section~\ref{sec:spectrum}, its spectrum is studied in detail
for various choices of su$_q$(2) Casimir operators and $q$ ranges.
The effect of the deformation on the corresponding wave functions
is determined in Section~\ref{sec:quadrupole} by calculating the
quadrupole moment in $l=0$ states. Finally,
Section~\ref{sec:conclusion} contains the conclusion. \par
%
%
\section{ su$_q$(2)-Invariant Schr\"odinger Equation}
\label{sec:equation}

Let
\begin{equation}
  H_q = -\frac{\hbar^2}{2\mu} \left(\frac{\partial^2}{\partial r^2} +
\frac{2}{r}
  \frac{\partial}{\partial r} - \frac{C_q}{r^2}\right) + \frac{1}{2} \mu
\omega^2 r^2
  \label{eq:H}
\end{equation}
be the Hamiltonian of a $q$-deformed three-dimensional harmonic oscillator in
spherical coordinates $r$, $\theta$, $\phi$. Here $C_q$ is the su$_q$(2) Casimir
operator, which we may take as
\begin{equation}
  C_q = J_+ J_- + \left[J_3 - \frac{1}{2}\right]_q^2 - \frac{1}{4},
\label{eq:C}
\end{equation}
where $[x]_q \equiv \left(q^x - q^{-x}\right)/\left(q - q^{-1}\right)$, and $q =
e^{w} \in $ R$^+$ or $q = e^{{\rm i}w} \in S^1$ (but different from a
root of
unity). The operators $J_3$, $J_+$, $J_-$, satisfying the su$_q$(2) commutation
relations
\begin{equation}
  \left[J_3, J_{\pm}\right] = \pm J_{\pm}, \qquad \left[J_+, J_-\right] =
  [2J_0]_q,
\end{equation}
are defined in terms of the angular variables by
\begin{eqnarray}
  J_3 & = & - {\rm i} \partial_{\phi}, \nonumber \\
  J_+ & = & - e^{{\rm i}\phi} \left( \tan(\theta/2) [T_1]_q\, q^{T_2}
          + \cot(\theta/2) q^{T_1}\, [T_2]_q\right), \nonumber \\
  J_- & = & e^{-{\rm i}\phi} \left( \cot(\theta/2) [T_1]_q\, q^{T_2}
          + \tan(\theta/2) q^{T_1}\, [T_2]_q\right),  \label{eq:J}
\end{eqnarray}
with $T_1 = - \frac{1}{2} \left(\sin \theta \partial_{\theta} - {\rm i}
\partial_{\phi}\right)$, $T_2 = - \frac{1}{2} \left(\sin \theta
\partial_{\theta} +
{\rm i} \partial_{\phi}\right)$~\cite{rideau,irac96}.\par
%
%
Instead of Equation~(\ref{eq:C}), we may alternatively use the operator
\begin{equation}
  C'_q = J_+ J_- + \left[J_3\right]_q \left[J_3-1\right]_q  \label{eq:C'}
\end{equation}
in Equation~(\ref{eq:H}), in which case the corresponding Hamiltonian will be
denoted by~$H'_q$.\par
%
%
The Hamiltonians $H_q$ and $H'_q$ remain invariant under su$_q$(2) since they
commute with $J_3$, $J_+$, $J_-$, and they coincide with the Hamiltonian of the
standard three-dimensional isotropic oscillator when $q=1$. For
simplicity's sake,
we shall henceforth adopt units wherein $\hbar = \mu = \omega = 1$.\par
%
%
In this representation, the stationary wave functions $\mid n l
m>_q$ can be written as
${\cal R}_{nlq}(r) Y_{lmq}(\theta,\phi)$,
%
%
 where $Y_{lmq}(\theta,\phi)$ are the $q$-spherical harmonics
 introduced in \cite{rideau}. We proved that they
  form an orthonormal set with respect to the scalar
  product~\cite{irac98}.

%
%
The radial wave functions are given by
\begin{equation}
  {\cal R}_{nlq}(r) = \left(\frac{2 (n!)}{\Gamma\left(\alpha_l + n +
  \frac{1}{2}\right)}\right)^{1/2} e^{- \frac{1}{2} r^2} r^{\alpha_l-1}
L_n^{\alpha_l -
  \frac{1}{2}}(r^2), \label{eq:S}
\end{equation}
where $L_n^{\alpha_l - \frac{1}{2}}\left(r^2\right)$ is an
associated Laguerre polynomial, and $\alpha_l
\in $ R$^ +$ is solution of the equation
\begin{equation}
  \alpha_l (\alpha_l  - 1) = C_q(l). \label{eq:alpha-equation}
\end{equation}
 The corresponding
energy
eigenvalues are
\begin{equation}
  E_{nlq} = 2n + \alpha_l + \frac{1}{2}, \qquad n = 0, 1, 2, \ldots.
\label{eq:E}
\end{equation}
The $\frac{1}{2}(N+1)(N+2)$-degeneracy of the isotropic oscillator energy
levels,
where $N = 2n+l$, is therefore lifted.\par
%
%
Similar results hold for the choice~(\ref{eq:C'}) for the Casimir
operator, the only change being the substitution of $C'_q(l) =
[l]_q [l+1]_q$ for $C_q(l)$. To distinguish the latter choice from
the former, we shall denote all quantities referring to it by
primed letters.\par
%
%
\section{ Spectrum of the su$_q$(2)-Invariant  Oscillator}
\label{sec:spectrum}
In the present section, we will study the condition $\alpha_l \in $ R$^+$ for the
existence of the radial wave functions~(\ref{eq:S}) and of the
corresponding energy
eigenvalues~(\ref{eq:E}), as well as the behaviour of the latter as
functions of $l$
and $q$ for the two choices (\ref{eq:C}), (\ref{eq:C'}) of Casimir
operators, and for
$q \in $ R$^+$ or $q \in S^1$.\par
%
%
$\clubsuit$ Let us first consider the case where $q =
e^{w}
\in $ R$^+$. Since the spectrum is clearly invariant under the
substitution $q
\to q^{-1}$, we may assume $q>1$, i.e., $w > 0$.\par
%
%
$\bullet$ In the $C_q$ case, for $l \ne 0$ only one root
$\alpha_{l+}$ of Equation~(\ref{eq:alpha-equation}) is positive and
therefore admissible, whereas for $l=0$, two roots $\alpha_{0+}$
and $\alpha_{0-}$ are admissible. Note that in the undeformed case
($q=1$), one root vanishes and has then to be discarded in
accordance with known results. For $q\ne1$, the spectrum therefore
comprises the energy eigenvalues
\begin{equation}\label{eq:l=0-C-R}
\begin{array}{lll}
 E_{n0q\pm} & = & 2n + 1 \pm \frac{1}{2 \cosh(w/2)},  \\
 &&\\
  E_{nlq} & = & 2n + 1 + \frac{\sinh((l+1/2)w)}{\sinh w}, \qquad l =
1, 2, \ldots.
\end{array}
\end{equation}
The appearance of an additional $l=0$ level was already observed in the
deformed Coulomb potential case~\cite{irac96}.\par
%
%
$\bullet$ In the $C'_q$ case, $\forall l$, only one root of
Equation~(\ref{eq:alpha-equation}) is positive and the spectrum
therefore comprises the same levels as in the undeformed case,
their energies being now
\begin{equation}
  E'_{nlq} = 2n + 1 + \frac{[4 \sinh(lw) \sinh((l+1)w) +
\sinh^2w]^{1/2}}
  {2\sinh w}, \qquad l = 0, 1, 2, \ldots.  \label{eq:l-C'-R}
\end{equation}
Note that the energy of the $l=0$ states, $E'_{n0q} = E_{n0} = 2n +
\frac{3}{2}$, is left undeformed  and, is therefore independent of $w$.
%
%
For $l\ne0$, $E_{nlq}$ and $E'_{nlq}$ are increasing functions of
$w$ in the neighbourhood of $w=0$, whereas for $l=0$, $E_{n0q+}$
and $E_{n0q-}$ have opposite behaviours. Moreover, for a given $w$
value, the influence of the deformation increases with $l$.
%

$\clubsuit$ Let us next consider the case where $q = e^{{\rm i}w}
\in S^1$. Owing to the invariance of the spectrum under the
substitution $q
\to q^{-1}$, we may now assume $0 < w < \pi$.\par
%
%
$\bullet$ For $C_q$, one finds two admissible roots if
$\gamma_q(l)\equiv 4
\sin^2w C_q(l) < 0$, but only one root $\alpha_{l+}$ if
$\gamma_q(l) \ge 0$. For instance, for $l=0$, there is a single
eigenvalue for any $w$ value,
\begin{equation}
  E_{n0q} = 2n + 1 + \frac{1}{2\cos(w/2)},  \label{eq:l=0-C-S}
\end{equation}
while for $l=1$, one obtains :
\begin{equation}
  \begin{array}{ll}
\mbox{two eigenvalues,}\
 E_{n1q\pm} = 2n + 1 \pm
\frac{4\cos^2(w/2) - 1}{2\cos(w/2)},
&\mbox{if\ } \frac{-7-\sqrt{17}}{16} < \cos w <
\frac{-7+\sqrt{17}}{16}, \\
&\\
 \mbox{one eigenvalue, }\
  E_{n1q} = 2n + 1 + \frac{4\cos^2(w/2) - 1}{2\cos(w/2)},
    & \mbox{if\ } \cos w \le \frac{-7-\sqrt{17}}{16}  \\
    &\mbox{or\ } \cos w \ge \frac{-7+\sqrt{17}}{16}.
  \end{array}
\end{equation}
\par
%
%
%
$\bullet$ For $C'_q$, one finds
that there are two admissible roots $\alpha'_{l+}$, $\alpha'_{l-}$ if $-
\sin^2w <
\gamma'_q(l) \equiv 4 \sin^2w C'_q(l) < 0$, only one $\alpha'_{l+}$
if either $\gamma'_q(l) \ge 0$ or $\gamma'_q(l) = - \sin^2w$ (in
which case $\alpha'_{l+} = 1/2$), or none if $\gamma'_q(l) < -
\sin^2w$. For instance, for $l=0$, one obtains a single
eigenvalue that coincides with the undeformed one, while for $l=1$,
one finds~:
\begin{equation}
  \begin{array}{ll}
 \mbox{two eigenvalues, }\ E'_{n1q\pm}
  =  2n + 1 \pm \frac{1}{2} \sqrt{1 + 8\cos w},
&\mbox{if }
           - \frac{1}{8} < \cos w < 0, \\
 \mbox{one eigenvalue, }\  E'_{n1q}  = 2n + 1 + \frac{1}{2} \sqrt{1 + 8\cos w},
 &\mbox{if }
           \cos w = - \frac{1}{8} \mbox{ or } \cos w \ge 0,\\
\mbox{none,}&\mbox{if}\cos w < -1/8.
           \end{array}
\end{equation}
%
%
In both the $C_q$ and $C'_q$ cases, for given $n$ and $l$ values
one always finds a single eigenvalue going into the undeformed one,
$E_{nl} = 2n + l + 3/2$, for $q
\to 1$. We can show that in a small enough
neighbourhood of $w=0$, for $l\ne0$, $E_{nlq}$ and $E'_{nlq}$ are
decreasing functions of $w$, whereas $E_{n0q}$ is increasing and
$E'_{n0q}$ remains constant.\par
%
%
For high $w$ values, the influence of the deformation on the
spectrum is rather striking and the levels get mixed in a very
complicated way.
%
%
In the $C'_q$ case, the situation is still more complex as some $l$
values may disappear on some intervals.
%
%
\section{ Quadrupole Moment in $l=0$ States}
\label{sec:quadrupole}

The purpose of the present section is to study the effect of the
deformation on the wave functions ${\cal R}_{nlq}(r)
Y_{lmq}(\theta,\phi)$ of the su$_q$(2)-invariant harmonic
oscillator, by determining the variation with $w$ of the
expectation value of some physically-relevant operator. For the
latter, we choose the electric quadrupole moment operator, and we
consider the quadrupole moment in the state $\mid n 0 0>_q$ which
is defined conventionally as
\begin{equation}
  Q_{n0q} = \langle n 0 0| (3z^2 - r^2) | n 0 0 \rangle_q.  \label{eq:Q}
\end{equation}
Equation~(\ref{eq:Q}) corresponds to the choice $C_q$ for the
su$_q$(2) Casimir operator. When using instead $C'_q$, the
quadrupole moment will be denoted by $Q'_{n0q}$.
%
%
As their undeformed counterpart vanishes, nonvanishing values of
$Q_{n0q}$ or $Q'_{n0q}$ will therefore be a direct measure of the
effect of the deformation. Note that since in the $C_q$ case, there
are two energy eigenvalues with $l=0$ and a given $n$ value for $q
\in$ R$^+$ (see Equation~(\ref{eq:l=0-C-R})), we have to
distinguish the corresponding quadrupole moments by a $\pm$
subscript.\par
%
%
 The quadrupole moment $Q_{n0q}$ can be factorized into radial and
angular matrix elements,
\begin{equation}
  Q_{n0q} = \langle n 0| r^2 |n 0\rangle_q \langle 00| \left(3 \cos^2\theta
- 1\right))|
  00\rangle_q.
\end{equation}
For any $l$ value, the former is simply obtained by replacing $l$ by
$\alpha_l - 1$
in the undeformed radial matrix element. Hence
\begin{equation}
  \langle n 0| r^2 |n 0\rangle_q = 2n + \alpha_0 + \frac{1}{2}.
\label{eq:def-radial}
\end{equation}
The calculation of the latter is more complicated as it implies the
use of the deformed angular scalar product~\cite{irac98},
\begin{eqnarray}
  \langle 00| (3 \cos^2\theta - 1)| 00\rangle_q & = & \frac{2
         \cosh^2w + 1}{\sinh^2w} - \frac{3 \cosh w}{w
\sinh w}, \qquad
         \mbox{if\ } q = e^{w} \in  \mbox{R}^+, \nonumber \\
  & = &  - \frac{2 \cos^2w + 1}{\sin^2w} + \frac{3 \cos w}{w
\sin w},
         \qquad \mbox{if\ } q = e^{{\rm i}w} \in S^1.
\label{eq:def-angular-bis}
\end{eqnarray}
The deformed quadrupole moment $Q_{n0q}$ is therefore an even function of
$w$, so that we may again restrict ourselves to $0 < w < \infty$ or
$0 < w <
\pi$ according to whether $q$ is real or complex.\par
%
%
%
{}For real $q$, it can be checked that $Q_{n0q-}$, $Q_{n0q+}$ and
$Q'_{n0q}$ are increasing positive functions of $w$. It should be
stressed that the undeformed quadrupole moments in $l\ne0$ states
and the deformed ones in $l=0$ states have opposite signs, and that
the effect of the choice of the Casimir and $\alpha_0$ root becomes
significant only for very large deformation.
%
%

{}For complex $q$, $Q_{n0q}$ and $Q'_{n0q}$ are decreasing negative
functions of $w$. Apart from the sign, which is now the same as
that of $Q_{nl}$ for $l\ne0$, the conclusions remain similar to
those for the real $q$ case.

%
%
\section{Conclusion}
\label{sec:conclusion}

In the present contribution, we did show that as those of the free
particle and of the Coulomb potential~\cite{irac96}, the
su$_q$(2)-invariant Schr\"odinger equation of the three-dimensional
harmonic oscillator can be easily solved not only for $q
\in $ R$^+$, but also for generic $q \in S^1$. It is worth stressing that we
have been working in the framework of the usual Schr\"odinger
equation (i.e., with no non-commuting objects contrary to some
other approaches~\cite{carow94}), but with wave functions belonging
to a Hilbert space different from the usual one, since the angular
part of the scalar product has been modified when going from su(2)
to su$_q$(2)~\cite{irac98}.\par
%
%
In the real $q$ case, we did show that the spectrum is rather similar to the
undeformed one, except that the energy levels are no more equidistant and that
their degeneracy is lifted. For a given $n$ value, the spacing between adjacent
levels corresponding to $l$ and $l+1$, respectively, increases with $l$ and with
the deformation. In addition, there appears a supplementary series of $l=0$
levels
when the Casimir operator $C_q$ is used. Apart from this, for small
deformations,
the results are rather insensitive to the choice made for the Casimir
operator.\par
%
%
In the complex $q$ case, we did show that the spectrum is more complicated
as for
$l\ne0$ and any $n$ value, there may exist 0, 1, or 2 levels according to the
deformation. The existence or inexistence of levels is also rather
sensitive to the
choice made for the Casimir operator. Close enough to $q=1$, there however
always
exists a single level going into the undeformed one for $q \to 1$. In that
region, the
spacing between adjacent levels corresponding to $l$ and $l+1$,
respectively, now
decreases with $l$ and with the deformation.\par
%
%
The closeness of our approach to the standard one did also allow us to study the
effect of the deformation on the wave functions and the Hilbert space spanned by
them. We did establish that it is rather strong as the quadrupole moment in the
$l=0$ states, which vanishes in the undeformed case, now assumes a positive
(resp.\ negative) value for $q \in $ R$^+$ (resp.\ $q \in S^1$) irrespective
of the
Casimir operator used.\par
%
%

%
%

\end{document}